\documentstyle[12pt]{article} \vsize=28.7truecm \hsize=23truecm
 \def\vbar{\mathchoice{\vrule height2.3ptdepth-.3ptwidth.12pt\kern-
 .10pt}
    {\vrule height6.3ptdepth-.3ptwidth.11pt\kern-.11pt}
    {\vrule height5.1ptdepth-.30ptwidth.8pt\kern-.8pt}
    {\vrule height4.1ptdepth-.24ptwidth.6pt\kern-.7pt}}
\def\reel{\hbox{I\hskip-2pt R}}

\def\<{\langle}
\def\>{\rangle}

\def\n{{\boldmath n}}

\def\<{\langle}
\def\>{\rangle}
\def\reel{\hbox{I\hskip -2pt R}}

\def\n{{\noindent}}
=0 in \headheight =0 cm

\newtheorem{theorem}{Theorem}[section]

\newtheorem{coro}[theorem]{Corollary}

\newtheorem{prop}[theorem]{Proposition}

\newenvironment{Pff}{{\bf Proof}}{\hfill $\Box$
\vspace*{0.2cm}}

\begin{document}
 \vspace{2cm}

\title{{\bf{Bayesian approach to cubic natural exponential families}}}
\author{Marwa Hamza and Abdelhamid Hassairi \footnote{Corresponding author.
 \textit{E-mail address: Abdelhamid.Hassairi@fss.rnu.tn}}
\\{\footnotesize{\it
Laboratory of Probability and Statistics. Sfax Faculty of
Sciences, B.P. 1171, Tunisia.}}}
\date{}
\maketitle \n $\overline{\hspace{15cm}}$\vskip0.3cm \n {\small {\bf
Abstract}} {\small For a natural exponential family (NEF), one can
associate in a natural way two standard families of conjugate
priors, one on the natural parameter and the other on the mean
parameter. These families of conjugate priors have been used to
establish some remarkable properties and characterization results of
the quadratic NEF's. In the present paper, we show that for a NEF,
we can associate a class of NEF's, and for each one of these NEF's,
we define a family of conjugate priors on the natural parameter and
a family of conjugate priors on the mean parameter which are
different of the standard ones. These families are then used to
extend to the Letac-Mora class of real cubic natural exponential
families the properties and characterization results related to the
Bayesian theory established for the quadratic natural exponential
families. }

\n {\small {\it{ Keywords:}} natural exponential family, variance
function, prior distribution, posterior expectation.\\
 $\overline{\hspace{15cm}}$\vskip1cm
\section{Introduction and preliminaries}
To make clear the motivations of the present paper, we first
recall some facts concerning the natural exponential families and
their variance functions. Our notations are the ones used by Letac
in [9]. Let $\mu$ be a positive radon measure on $\reel$, and
denote by
\begin{equation}\label{MI1}
L_\mu(\lambda)=\displaystyle\int \exp(\lambda x)\mu(dx)
\end{equation}
its Laplace transform. Let ${\mathcal{M}}(\reel)$ be the set of
measures $\mu$ such that the set
\begin{equation}
\Theta(\mu)=\textrm{interior}\{\lambda \in\reel; \ L_\mu
(\lambda)<+\infty\}
\end{equation}
is non empty and $\mu$ is not Dirac measure. The cumulant function
of an element $\mu$ of ${\mathcal{M}}(\reel)$ is the function
defined for $\lambda\in\Theta(\mu)$ by
$$ k_\mu (\lambda)=\ln L_\mu (\lambda).$$ \vskip0.1cm\n To
each $\mu$ in ${\mathcal{M}}(\reel)$ and $\lambda$ in
$\Theta(\mu)$, we associate the probability distribution on
$\reel$
\begin{equation}\label{a1} P(\lambda,\mu)(dx)=\exp(\lambda x-k_\mu
(\lambda))\mu(dx).\end{equation}
 The set\vskip0.1cm$\hfill
 F(\mu)=\{P(\lambda,\mu); \ \lambda\in\Theta(\mu
 )\}\hfill$
\vskip0.3cm\n is called the natural exponential family (NEF)
generated by $\mu$.\vskip0.2cm\n The function $k_\mu$ is strictly
convex and analytic. Its first derivative $k'_\mu$ defines a
diffeomorphism between $\Theta(\mu)$ and its image $M_{F(\mu)}$.
Since $k'_\mu (\lambda)=\displaystyle\int x\ P(\lambda,\mu) (dx)$,
$M_{F(\mu)}$ is called the domain of the means of $F(\mu)$. The
inverse function of $k'_\mu$ is denoted by $\psi_\mu$ and setting
\begin{equation}\label{b1} P(m,F(\mu))=P(\psi_\mu(m),\mu)\end{equation} the probability of $F(\mu)$ with mean
$m$, we have $$F(\mu)=\left\{P(m,F(\mu)); \ m\in
M_{F(\mu)}\right\},$$ which is the parametrization of $F(\mu)$ by
the mean.\\ The variance of $P(m,F(\mu))$ is denoted by $V_{F(\mu)}
(m)$ and the map defined from $M_{F(\mu)}$ into $L_s(\reel)$, the
set of symmetric function in $\reel$, by
$$ m\longmapsto V_{F(\mu)} (m)=k''_\mu(\psi_\mu(m))=(\psi'_\mu(m))^{-1}$$
is called the variance function of the NEF $F(\mu)$ generated by
$\mu$. We also say that $\mu$ is a basis of $F(\mu)$. An important
feature of $V_{F(\mu)}$ is that it characterizes the natural
exponential family $F(\mu)$ in the following sense: If $F(\mu)$
and $F(\nu)$ are two NEF's such that $V_{F(\mu)} (m)$ and
$V_{F(\nu)}(m)$ coincide on
a nonempty open subset of $M_{F(\mu)} \cap M_{F(\nu)},$ then $F(\mu)=F(\nu)$.\\
Now, for $\mu\in\mathcal{M}(\reel)$ the J{\o}rgensen set of $\mu$ or
of $F(\mu)$ is defined by
$$\Lambda(\mu)=\{\lambda>0;\ \exists\
\mu_{\lambda} ;\
L_{\mu_{\lambda}}(\theta)=(L_{\mu}(\theta))^{\lambda}\ and\
\Theta(\mu_{\lambda})=\Theta(\mu) \}.$$

\n $\Lambda(\mu)$ is stable under addition which means that if
$\lambda$ and $\lambda'$ are in $\Lambda(\mu)$, then
$\lambda+\lambda'$ are in $\Lambda(\mu)$, and
$\mu_{\lambda+\lambda'}=\mu_{\lambda}*\mu_{\lambda'}.$

\n For all $\lambda$ in $\Lambda(\mu)$ we have
$$M_{F(\mu_{\lambda})}=\lambda M_{F(\mu)}\ \ and\ \
V_{F(\mu_{\lambda})}(m)=\lambda V_{F(\mu)}(\frac{m}{\lambda}).$$
Several classifications of NEFs based on the form of the variance
function have been realized in the last three decades. The most
interesting classes of real NEF's are the class of quadratic NEFs,
i.e., the class of NEF's such that the variance function is a
polynomial in the mean of degree less than or equal to 2
characterized by Morris \cite{Morris(1982)}, and the class of cubic
NEF's, i.e., the class of NEF's such that the variance function is a
polynomial in the mean of degree less than or equal to 3
characterized by Letac and Mora \cite{Letac and Mora (1990)}. Recall
that up to affine transformations and power of convolution the class
of quadratic NEF's contains six families: the gaussian, the Poisson,
the gamma, the binomial, the negative binomial, an the hyperbolic
family. The class of cubic NEF's contains, besides the quadratic
ones, six other families, the most famous is the inverse Gaussian
distribution with variance function $V(m)=m^{3}$. It is worth
mentioning here that multivariate versions of these classes have
also been defined and completely described. For instance, Casalis
\cite{Casalis(1996)} has described the so-called class of
multivariate simple quadratic NEFs and Hassairi \cite{H} has
described the class of multivariate simple cubic NEFs which are
respectively the generalizations of the real quadratic and real
cubic NEF's. The fact that the variance function of a family is
quadratic or cubic is not only a question of form but it corresponds
to some interesting analytical characteristic properties. Indeed,
the Morris class of quadratic NEF's has some characterizations
involving orthogonal polynomials due to Fiensilver \cite{Feinsilver
(1986)}. These characterizations have been extended to the Letac and
Mora class of real cubic NEF's by Hassairi and Zarai \cite{Hassairi
and Zarai(2004)} using a notion of 2-orthogonality of a sequence of
polynomials. Other remarkable characterizations of the quadratic
NEF's are related to the Bayesian theory. For instance,  given a NEF
$F(\mu)$, Diaconis and Ylvisaker $\cite{Diconis and
Ylvisaker(1979)}$ have considered the standard family $\Pi$ of
priors on the natural parameter $\lambda$ defined by
\begin{equation}\label{1}
\pi_{t_{1},m_{1}}(d\lambda)=C_{t_{1},m_{1}}\ \exp(t_{1}
m_{1}\lambda-t_{1}k_{\mu}(\lambda))\
\mathbf{1}_{\Theta(\mu)}(\lambda)\ \textit{d}\lambda\end{equation}
where $t_{1}>0$, $m_{1}$ is in $M_{F(\mu)}$, and $C_{t_{1},m_{1}}$
is a normalizing constant. This distribution is in fact a particular
case of the so called implicit distribution on the parameter of a
statistical model introduced in \cite{H.M.K}. They have shown that
if $X$ is a random variable distributed according to
$P(\lambda,\mu)$ (see (\ref{a1})), then the only conjugate family of
prior distributions on $\lambda$ that gives a linear posterior
expectation of $k'_{\mu}(\lambda)$ given $X$ is the standard one
$\Pi$ (see also \cite{CR}). Consonni and Veronese \cite{CV} have
considered another family $\widetilde{\Pi}$ of prior distributions
$\widetilde{\pi}_{t_{1},m_{1}}$ on the mean parameter $m$ defined
also for $t_{1}>0$ and $m_{1}$ in $M_{F(\mu)}$ by
\begin{equation}\label{a2}\widetilde{\pi}_{t_{1},m_{1}}(dm)=\widetilde{C}_{t_{1},m_{1}}\exp(t_{1}m_{1}\psi_{\mu}(m)-t_{1}k_{\mu}(\psi_{\mu}(m)))\
\mathbf{1}_{\textit{M}_{\textit{F}(\mu)}}(\textit{m})\textit{dm}.\end{equation}
They have shown that the fact that $\widetilde{\Pi}$ contains
$k'_{\mu}(\Pi)$ characterizes the quadratic NEFs. These authors have
also shown that if the prior on the mean parameter $m$ is
$\widetilde{\pi}_{t_{1},m_{1}} $, then under some conditions on the
support of $\mu$, the NEF $F(\mu)$ is quadratic if and only if the
posterior expectation of $k'_{\mu}(\lambda)$ is a linear function of
the sample mean. We also mention that Diaconis and Yilvisaker
$\cite{Diconis and Ylvisaker(1979)}$ have shown that if the standard
prior on $\lambda$ is given by $\pi_{t_{1},m_{1}}$ with $t_{1}>0$
and $m_{1}$ is in $M_{F(\mu)}$, then the expectation of
$k'_{\mu}(\lambda)$ is equal to $m_{1}$, that is
\begin{equation} \label{a18} C_{t_{1},m_{1}}\displaystyle\int
k'_{\mu}(\lambda)\ \exp(t_{1}m_{1}\lambda-t_{1}k_{\mu}(\lambda))\
\mathbf{1}_{\Theta(\mu)}(\lambda)\textit{d}\lambda =
\textit{m}_{1},
\end{equation} or equivalently in terms of the mean parameter \begin{equation}\label{a19}C_{t_{1},m_{1}}\displaystyle\int
\displaystyle\frac{m}{V_{F(\mu)}(m)}\exp(t_{1}m_{1}\psi_{\mu}(m)-t_{1}k_{\mu}(\psi_{\mu}(m)))\
\mathbf{1}_{\textit{M}_{\textit{F}(\mu)}
}(\textit{m})\textit{dm}=\textit{m}_1.\end{equation} A natural
question within this approach is if one can extend the properties
and characterization results concerning the quadratic NEF's and
related to the Bayesian theory to the Letac-Mora class of real cubic
NEFs. The aim of the present paper is to give an answer to this
question. We first introduce, for a given NEF $F(\nu)$ and $\beta$
in some interval of $\reel$ containing 0, a NEF $F^{\beta}(\nu)$
such that $F^{0}(\nu)=F(\nu)$. We then define a family $\Pi^{\beta}$
of prior distributions on the natural parameter $\theta$ and a
family $\widetilde{\Pi}^{\beta}$ of prior distributions on the mean
parameter $m$ which may be seen as generalizations of the families
$\Pi$ and $\widetilde{\Pi}$ defined above, since $ \Pi=\Pi^{0}$ and
$\widetilde{\Pi}=\widetilde{\Pi}^{0}$. After proving that for each
$\beta$, the family $\widetilde{\Pi}^{\beta}$ is a conjugate family
of prior distributions with respect to the NEF $F^{\beta}(\nu)$, we
show that a cubic NEF $F(\nu)$ is characterized by the fact that
there exists a $\beta$ such that the posterior expectation of
$\displaystyle\frac{k'_{\nu}(\theta)}{1-\beta k'_{\nu}(\theta)}$ is
linear when the prior on $\theta$ is $\pi^{\beta}_{t,m_{0}}.$ We
also show that a cubic NEF $F(\nu)$ is characterized by a
differential equation verified by the cumulant function $k_{\nu}$. A
third characterization of a real cubic NEF is realized when the
family of priors $\widetilde{\Pi}^{\beta}$ contains the family
$k'_{\nu}(\Pi^{\beta})$.  The restriction of all these results to
the subclass of quadratic NEF's leads to the results of Diconis and
Ylvisaker\cite{Diconis and Ylvisaker(1979)} and Consonni and
Veronese\cite{CV}. The results of the paper are illustrated by an
example.
\section{Main results}
In this section, we state and prove our main results. Our
considerations will be restricted to regular NEFs, so that the
domain of the means of a NEF is equal to the interior of the convex
hull of its support. This property of regularity is satisfied by all
the most common NEF's. An important fact which will be crucial in
our proofs is that up to affine transformations and powers of
convolution, a cubic natural exponential family may be obtained from
a quadratic one by the so-called action of the linear group
$GL(\reel^{2})$ on the real families. Originally, the action of the
linear group includes the affine transformations and powers of
convolution, however since these transformations preserve the class
of quadratic NEF's and the class of cubic NEF's, we will focus on
the facts which we need here, for more precise statements in this
connection, refer to Hassairi\cite{H}. Let $F(\nu)$ and $F(\mu)$ be
two real NEFS's. Suppose that there exists a $\beta$ in $\reel$ such
that the set
$$(M_{F(\nu)})_{\beta}=\{m\in M_{F(\nu)}; 1+\beta m>0\}$$
is nonempty and for $m$ in $(M_{F})_{\beta}$,
\begin{equation}\label{v}V_{F(\nu)}(m)=(1+\beta m)^{3}\ V_{F(\mu)}\left(\displaystyle\frac{m}{1+\beta m}\right)\
,\end{equation} then we write $F(\nu)=T_{\beta}(F(\mu))$. This
defines an action on the natural exponential families, so that we
have
 $T_{\beta}T_{\beta'}=T_{\beta+\beta'}$ and
$F(\nu)=T_{\beta}(F(\mu))$ is equivalent to
$F(\mu)=T_{-\beta}(F(\nu))$. We also mention that
$F(\nu)=T_{\beta}(F(\mu))$ may be expressed in terms of the
cumulant functions of the generating measures by\begin{equation} \label{e6}\left\{\begin{array}{ccc}k_{\nu}(\theta)&=&k_{\mu}(\lambda)\\\
\\ \theta&=&-\beta
k_{\mu}(\lambda)+\lambda\end{array}\ \ \right.\end{equation} or
equivalently by

\n \

\begin{equation}\label{ee} \left\{\begin{array}{ccc}
k_{\mu}(\lambda)&=&k_{\nu}(\theta)\\\
\\\lambda&=&\beta
k_{\nu}(\theta)+\theta\end{array}\ \ \right.
\end{equation}
An important fact is that the relation (\ref{ee}) between the
cumulent functions may be explicitly given it terms of the measures
themselves. In fact if the $\alpha$-power of convolution
$\nu_{\alpha}$ of $\nu$ is written as
$\nu_{\alpha}(dx)=h(\alpha,x)\sigma(dx)$, where $\sigma(dx)$ is
either Lebesgue measure or a counting measure, then the measure
$$\mu(dx)=\frac{1}{1-\beta x}h(1-\beta x,x)\mathbf{1}_
{\Lambda(\nu)}(\textit{1}-\beta \textit{x})\sigma(\textit{dx})$$
satisfies (\ref{e6}) and generates the family
$T_{-\beta}(F(\nu)).$ This measure $\mu$ will be denoted
$T_{-\beta}(\nu)$ and the family $F(\mu)=T_{-\beta}(F(\nu))$ will
be denoted
$F^{\beta}$. \\
\n We mention here that if $F(\nu)$ is a cubic NEF, there exists
$\beta$ in $B_{F(\nu)}$ and a quadratic NEF $F(\mu)$ such that
$F(\nu)=F(T_{\beta}(\mu))$ or
equivalently, $F(\mu)=F(T_{-\beta}(\nu)).$\\
Besides the restriction to half lines for the domain of the means in
the definition of $(M_{F(\nu)})_{\beta}$, we also define for
$\nu\in\mathcal{M}(\reel)$ and $\beta\in\reel$, the sets
$$H_{\beta}(\nu)=\{x\in\reel;\  1+\beta x\in \Lambda(\nu)\},$$ and
$$B_{F(\nu)}=\{\beta\in\reel;\ \nu(H_{\beta})>0\}.$$ We have
the following preliminary result.
\begin{prop} Let $F(\nu)$ be a regular NEF. Then $$\beta \in B_{F(\nu)} \ \ if\ \ and\ \ only \ \ if \ \
(M_{F(\nu)})_{\beta}\neq\emptyset.$$
\end{prop}
\begin{Pff}\ \ We will make a reasoning for $\beta\geq0$, the case
$\beta<0$ may be done in a similar way. Suppose that there exists
$m_{0}$ in $(M_{F(\nu)})_{\beta}$, that is $m_{0}\in M_{F(\nu)}$ and
$1+\beta m_{0}>0$. As $M_{F(\nu)}$ is equal to the interior of the
convex hull of $supp(\nu)$, there exist $x_{0}$ in $supp(\nu)$ such
that $x_{0}>m_{0}.$ This with the fact that $1+\beta m_{0}>0$ imply
that $1+\beta x_{0}>0$. Thus $H_{\beta}$ is an open set which
contains an element of $supp(\mu)$. It follows that
$\nu(H_{\beta})>0$ and $\beta$ is in $B_{F(\nu)}.$
\\Conversely, if $\beta$ is in $B_{F(\nu)}$, then $\nu(H_{\beta})>0$.
Since $H_{\beta}$ is an open set, this implies that it contains an
element $x_{0}$ of $supp(\nu)$. We have on the one hand that
$1+\beta x_{0}>0$ so that there exists $\varepsilon>0$ such that
$1+\beta x_{0}-\beta\varepsilon>0$. On the other hand, as
$M_{F(\nu)}$ is equal to the interior of the convex hull of
$supp(\nu)$, there exists $m_{0}$ in $(M_{F(\nu)})$ such that
$|m_{0}-x_{0}|<\varepsilon$. From this we deduce that $m_{0}$ is in
$(M_{F(\nu)})_{\beta}$.\end{Pff}\\
We now use the natural parametrization and the parametrization by
the mean of the original family $F(\nu)$ to give two
parameterizations of the family $F(\mu)$. These parameterizations
are, for $\beta\neq 0$, different of the usual parameterizations of
$F(\mu)$. In fact, for $\theta\in\Theta(\nu)$, we write
$$P(\beta,\theta,\nu)(dx)=\exp\{(\theta+\beta k_{\nu}(\theta)) x-k_{\nu}(\theta)\}\
T_{-\beta}(\nu)(dx).$$ Similarly, parameterizing by $m\in
M_{F(\nu)}$, we write
$$P(\beta,m,F(\nu))(dx)=\exp\{(\psi_{\nu}(m)+\beta k_{\nu}(\psi_{\nu}(m))) x-k_{\nu}(\psi_{\nu}(m))\}\
T_{-\beta}(\nu)(dx).$$ Thus we have that
$$F(\mu)=F^{\beta}=\{P(\beta,\theta,\nu)(dx);\ \theta\in\Theta(\nu)\}=\{P(\beta,m,F(\nu))(dx);\ m\in M_{F(\nu)}\}.$$
Accordingly, we define for $\beta$ in $B_{F(\nu)}$ two families of
prior distributions. Let
$$(\Theta)_{\beta}=\{\theta\in\Theta(\nu);\ 1+\beta
k'_\nu(\theta)>0\}.$$ Then we have that
$(M_{F(\nu)})_{\beta}=k'_\nu((\Theta)_{\beta}),$ and we define for
$t>0\ and\ m_{0}\in (M_{F(\nu)})_{\beta}$,
\begin{equation}\label{d}\pi^{\beta}_{t,m_{0}}(\textit{d}\theta)=C^{\beta}_{t,m_{0}}\ (1+\beta
k'_\nu(\theta))\ \exp(tm_{0}\theta-tk_{\nu}(\theta))\
\mathbf{1}_{(\Theta)_{\beta}}(\theta)
\textit{d}\theta,\end{equation} and
$$\Pi^{\beta}=\{\pi^{\beta}_{t,m_{0}};\ t>0\ and\ m_{0}\in
(M_{F(\nu)})_{\beta}\},$$ This family comes in fact from the
standard family $\Pi$ defined in (\ref{1}) using (\ref{ee}). The
normalizing constant $C^{\beta}_{t,m_{0}}$ is then well defined for
$t>0\ and\ m_{0}\in (M_{F(\nu)})_{\beta}$.\\Besides this family of
priors on the parameter $\theta$, we define a family of priors on
the parameter $m$. Always for $t>0\ and\ m_{0}\in
(M_{F(\nu)})_{\beta}$, we consider the probability distribution
$$\widetilde{\pi}^{\beta}_{t,m_{0}}(dm)=\widetilde{C}^{\beta}_{t,m_{0}}\ (1+\beta
m)^{-2}\ \exp(tm_{0}\psi_{\nu}(m)-tk_{\nu}(\psi_{\nu}(m)))\
\mathbf{1}_{(\textit{M}_{F(\nu)})_{\beta}}(\textit{m})\textit{dm},$$
where $\widetilde{C}^{\beta}_{t,m_{0}}$ is a normalizing constant.
It is the image of $\widetilde{\pi}_{t_{1},m_{1}}$ defined in
(\ref{a2}) by the map $m'\longmapsto
\displaystyle\frac{m'}{1-\beta m'}.$ The family of priors on $m$
is then
$$\widetilde{\Pi}^{\beta}=\{\widetilde{\pi}^{\beta}_{t,m_{0}};\ t>0\ and\ m_{0}\in
(M_{F(\nu)})_{\beta}\},$$
 \n Next, we prove that these families are conjugate families of
 prior distributions.
\begin{prop}\label{conj}i) The family $\Pi^{\beta}$ is conjugate family of prior
distributions with respect to the NEF $F^{\beta}$ parameterized by
the natural parameter $\theta$.

ii) The family $\widetilde{\Pi}^{\beta}$ is a conjugate family of
prior distributions with respect to the NEF $F^{\beta}$
parameterized by the mean parameter $m$.
\end{prop}
\begin{Pff} i) Suppose that
$X$ is a random variable with distribution
$P(\beta,\theta,\nu)(dx)$ and that $\pi^{\beta}_{t,m_{0}}$ is a
prior distribution on the parameter $\theta$. Then the posterior
distribution is
$$\displaystyle\frac{C^{\beta}_{t,m_{0}}\ (1+\beta k'_{\nu}(\theta))\
\exp((tm_{0}+x)\theta-(t+1-\beta x)k_{\mu}(\theta))\
\mathbf{1}_{(\Theta)_{\beta}}(\theta)\ }{\displaystyle\int
C^{\beta}_{t,m_{0}}\ (1+\beta k'_{\nu}(\theta))\
\exp((tm_{0}+x)\theta-(t+1-\beta x)k_{\mu}(\theta))\
\mathbf{1}_{(\Theta)_{\beta}}(\theta)\ \textit{d}\theta\ }.$$ If
we set $t_{2}=t+1-\beta x $ and
$m_{2}=\displaystyle\frac{tm_{0}+x}{t+1-\beta x}$, then this
distribution is noting but $\pi^{\beta}_{t+1-\beta
X,(tm_{0}+X)/(t+1-\beta X)},$ and it belongs to $\Pi^{\beta}$. In
fact, since $t>0$, $T_{-\beta}(\nu)$ is concentrated on $\{1-\beta
x>0\}$, and $m_{0}$ is in $(M_{F(\nu)})_{\beta}$, we have that
$t_{2}>0$ and $ 1+\beta m_{2}=\displaystyle\frac{1+t(1+\beta
m_{0})}{t+1-\beta X}>0,$ so that $m_{2}$ is in $(M_{F})_{\beta}$.\\
ii) Suppose now that $Y$ is a random variable $P(\beta,m,F(\nu))$
distributed and that the prior on the mean parameter $m$ is
$\widetilde{\pi}^{\beta}_{t,m_{0}}$, then the posterior distribution
of $m$ is given by
$$\displaystyle\frac{\widetilde{C}^{\beta}_{t,m_{0}}\ (1+\beta m)^{-2}\
\exp((tm_{0}+y)\psi_{\mu}(m)-(t+1-\beta y)k_{\mu}(\psi_{\mu}(m)))\
\mathbf{1}_{(\textit{M}_{F(\nu)})_{\beta}}(\textit{m})
}{\displaystyle\int \widetilde{C}^{\beta}_{t,m_{0}}\ (1+\beta
m)^{-2}\ \exp((tm_{0}+y)\psi_{\mu}(m)-(t+1-\beta
y)k_{\mu}(\psi_{\mu}(m)))\
\mathbf{1}_{(\textit{M}_{F(\nu)})_{\beta}}(\textit{m})\
\textit{dm}}.$$ This with $t_{3}=t+1-\beta y$ and
$m_{3}=\displaystyle\frac{tm_{0}+y}{t+1-\beta y}$, is equal to
$\widetilde{\pi}^{\beta}_{t+1-\beta Y,(tm_{0}+Y)/(t+1-\beta Y)}$,
with the required conditions to belong to
$\widetilde{\Pi}^{\beta}$.
\end{Pff} \begin{coro}\label{coro1} Let $(X_{1},...,X_{n})$ be a sample
$P(\beta,\theta,\nu)$-distributed and consider
$\pi^{\beta}_{t,m_{0}}$ as the prior on $\theta$. Then the
posterior distribution of $\theta$ given $X_{1},...,X_{n}$ is
$$\pi^{\beta}_{t+n-\beta n
\overline{X},(tm_{0}+n\overline{X})/(t+n-\beta n \overline{X})}.
$$
\end{coro}
\begin{Pff}\ It is easy to see that the distribution of the random
vector $(\theta, X_{1},...,X_{n})$ is
$$C^{\beta}_{t,m_{0}}\
(1+\beta k'_{\nu}(\theta))\
\exp((tm_{0}+n\overline{x})\theta-(t+n(1-\beta
\overline{x}))k_{\nu}(\theta))\
\mathbf{1}_{(\Theta)_{\beta}}(\theta)\
\prod_{\textit{i}=\textit{1}}^{\textit{n}}\textit{T}_{-\beta}(\nu)(\textit{dx}_{\textit{i}})
\textit{d}\theta.$$With the same technic used in Proposition
\ref{conj}, we deduce that the posterior distribution of $\theta$
given $X_{1},...,X_{n}$ is
$$\pi^{\beta}_{t+n-\beta n\overline{X},(tm_{0}+n\overline{X})/(t+n-\beta n\overline{X})}.$$\end{Pff}
\begin{prop}\label{prop}Let $F(\nu)$ be a cubic NEF. Then there exists $\beta$ in $B_{F(\nu)}$ such
that, if the prior on $\theta$ is $\pi^{\beta}_{t,m_{0}}$, then
$$E \left(
  \frac{k'_{\nu}(\theta)}{1+\beta k'_{\nu}(\theta)}
\right)= \frac{m_{0}}{1+\beta m_{0}}.$$\end{prop}
\begin{Pff}
Since $F(\nu)$ is cubic then there exist $\beta$ in $B_{F(\nu)}$
and a quadratic NEF $F(\mu)$ such that $F(\nu)=T_{\beta}(F(\mu))$
 and $\nu=T_{\beta}(\mu).$ Using (\ref{ee}) it is easy to see that if
the prior on $\theta$ in $(\Theta)_{\beta}$ is
$\pi^{\beta}_{t,m_{0}}$ then the prior of $\lambda$ in $\Theta(\mu)$
is the standard $\pi_{t_{1},m_{1}}$ with $t_{1}=t(1+\beta m_{0})$
and $m_{1}=\displaystyle\frac{m_{0}}{1+\beta m_{0}}.$ Moreover we
have
$C^{\beta}_{t,m_{0}}=C_{t(1+\beta m_{0}),m_{0}/(1+\beta m_{0})}.$\\
It follow that

\begin{eqnarray*}E\left(\frac{k'_{\nu}(\theta)}{1+\beta
k'_{\nu}(\theta)}\right)&=&\displaystyle\int C^{\beta}_{t,m_{0}}\
\frac{k'_{\nu}(\theta)}{1+\beta k'_{\nu}(\theta)}\ (1+\beta
k'_{\nu}(\theta))\ \exp(tm_{0}\theta-tk_{\nu}(\theta))\
\mathbf{1}_{(\Theta)_{\beta}}(\theta)\textit{d}\theta\\
&=&\displaystyle\int C^{\beta}_{t,m_{0}}\ k'_{\mu}(\lambda)\
\exp(tm_{0}(\lambda-\beta k_{\mu}(\lambda))-tk_{\mu}(\lambda))\
\mathbf{1}_{\Theta(\mu)}(\lambda)\textit{d}\lambda\\
&=& \displaystyle\int\ C_{t_{1},m_{1}}k'_{\mu}(\lambda)\
\exp(t_{1}m_{1}\lambda-t_{1}k_{\mu}(\lambda))\
\mathbf{1}_{\Theta(\mu)}(\lambda)\textit{d}\lambda.
\end{eqnarray*} Invoking
(\ref{a18}) we get
$$E \left(\frac{k'_{\nu}(\theta)}{1+\beta
k'_{\nu}(\theta)}\right)=\frac{m_{0}}{1+\beta m_{0}}.$$
\end{Pff}

\n Now we give a characterization of the cubic NEFs which is based
on the linearity of the posterior expectation.
\begin{theorem} Let $\nu$ be in $M(\reel)$.
\begin{enumerate}\item If $F(\nu)$ is cubic then there exists $\beta$ in
$B_{F(\nu)}$ such that for all $n\geq 1$, if $X_{1},...,X_{n}$ is
a sample with distribution $P(\beta,\theta,\nu)$ and the prior on
the natural parameter $\theta$ is $\pi^{\beta}_{t,m_{0}}$, then
$E\left(\frac{k'_{\nu}(\theta)}{1+\beta k'_{\nu}(\theta)}\mid
X_{1},...,X_{n}\right)$ is linear.\\ \item The converse is true if
we assume that $supp(T_{-\beta}(\nu))$ contains an open interval
in $\reel$ or is a finite or denumerably infinite subset of
$]-\infty,0[$, alternatively $]0,+\infty[$, with
$T_{-\beta}(\nu)\{0\}>0$.\end{enumerate}\end{theorem}
\begin{Pff}

\n \begin{enumerate}\item From Proposition \ref{prop} and Corollary
\ref{coro1}, we easily get
$$E\left(\frac{k'_{\nu}(\theta)}{1+\beta k'_{\nu}(\theta)}\mid
X_{1},...,X_{n}\right)=\frac{tm_{0}+n\overline{X}}{t(1+\beta
m_{0})+n} ,$$which is linear in $\overline{X}.$\\ \item
Conversely, suppose that there exists $\beta$ in $B_{F(\nu)}$ such
that $E\left(\frac{k'_{\nu}(\theta)}{1+\beta k'_{\nu}(\theta)}\mid
X_{1},...,X_{n}\right)$ is linear for all $n$ and a sample
$X_{1},...,X_{n}$ with distribution $P(\beta,\theta,\nu)(dx)$.
Consider $\mu=T_{-\beta}(\nu)$. Then using (\ref{ee}), we have
that the distribution $P(\beta,\theta,\nu)(dx)$ of
$X_{1},...,X_{n}$ is equal to $P(\lambda,\mu)(dx)$ which is
nothing but an other parametrization involving $\mu.$ We also have
that
$$E\left(k'_{\mu}(\lambda)\mid X_{1},...,X_{n}\right)=E\left(\frac{k'_{\nu}(\theta)}{1+\beta k'_{\nu}(\theta)}\mid
X_{1},...,X_{n}\right)$$ which is linear in $X_{1},...,X_{n}$,
from the hypothesis.
\\On the other hand, as the prior on the natural parameter $\theta$
is assumed to be $\pi^{\beta}_{t,m_{0}}$, we get as prior on
$\lambda\in \Theta(\mu)$ the standard prior distribution given by
$$\pi_{t_{1},m_{1}}(\textit{d}\lambda)=
C_{t_{1},m_{1}}\ \exp(t_{1}m_{1}\lambda-t_{1}k_{\mu}(\lambda))\
\mathbf{1}_{\Theta(\mu)}(\lambda)(\textit{d}\lambda),$$ with
$t_{1}=t(1+\beta m_{0})>0$,
$m_{1}=\displaystyle\frac{m_{0}}{1+\beta
m_{0}}\in M_{F(\mu)}$ and $C_{t_{1},m_{1}}=C^{\beta}_{t,m_{0}}.$\\
As we have that $supp(\mu)=supp(T_{-\beta}(\nu))\subset
supp(\nu)\cap\{1-\beta x\in \Lambda(\nu)\}.$ the assumptions on
$supp(T_{-\beta}(\nu))$ imply that $supp(\mu)$ satisfies
hypotheses (H1) or (H2) of Theorem 1.1 of Consonni and Veronese.
According to this and to the linearity of the conditional
expectation of the mean parameter of $F(\mu)$, we deduce that this
NEF is a quadratic. It follows that $F(\nu)=T_{\beta}(F(\mu))$ is
a cubic NEF.
\end{enumerate}
\end{Pff}\\
In the following theorem, we give a second characterization of the
Letac-Mora class of real cubic NEFs.
\begin{theorem}\label{T1}
Let $\nu$ be in $\mathcal{M}(\reel)$. $F(\nu)$ is cubic if and
 only if there exist $\beta$ in $B_{F(\nu)}$ and $(a,b,c)\in\reel^{3}$
such that for all $m$ in $M_{F(\nu)},$
\begin{equation}\label{e1}V_{F(\nu)}(m)=(1+\beta \textit{m})^{3}\
\exp{(a\psi_{\nu}(m)+bk_{\nu}(\psi_{\nu}(m))+c)}.\end{equation}\end{theorem}
Note that (\ref{e1}) may be expressed in terms of the cumulant
function as there exist $\beta$ and $(a,b,c)\in\reel^{3}$ such
that for all $\theta$ in $\Theta_{\beta}$
$$k''_{\nu}(\theta)=(1+\beta k'_{\nu}(\theta))^{3}\
\exp{(a\theta+bk_{\nu}(\theta)+c)},$$ that is the cumulant
function is solution of some Monge-Amp\`{e}re equation
(see\cite{Monge-Ampere}).
\begin{Pff}
Suppose that $F(\nu)$ is cubic, then there exist $\beta$ in
$B_{F(\nu)}$ and a quadratic NEF $F(\mu)$ such that
$F(\nu)=T_{\beta}(F(\mu))$ or equivalently
$F(\mu)=T_{-\beta}(F(\nu))$. Then it follows from (\ref{v}) that
$$V_{F(\nu)}(m)=(1+\beta m)^{3}\ V_{F(\mu)}\left(\frac{m}{1+\beta
m}\right).$$ It is known (see \cite{Casalis(1996)}) that for the
quadratic NEF $F(\mu)$ there exists $(a',b',c')\in\reel^{3}$ such
that for all $m'$ in $M_{F(\mu)}$
$$V_{F(\mu)}(m')=\exp(a'\psi_{\mu}(m')+b'k_{\mu}(\psi_{\mu}(m'))+c').$$
Writing (\ref{ee}) in terms of the mean parameters we get
\begin{equation}\label{psi}\left\{\begin{array}{ccc}
k_{\mu}\left(\psi_{\mu}(\frac{m}{1+\beta m})\right)&=&k_{\nu}(\psi_{\nu}(m))\\\
\\\psi_{\mu}\left(\frac{m}{1+\beta m}\right)&=&\beta
k_{\nu}(\psi_{\nu}(m))+\psi_{\nu}(m)\end{array}\ \
\right.\end{equation}Therefore
$$V_{F(\nu)}(m)=(1+\beta \textit{m})^{3}\
\exp(a\psi_{\nu}(m)+b k_{\nu}(\psi_{\nu}(m))+c),$$ with
$\label{aa1}a=a',\ \ b=b'+\beta a'\ \ and\ \ c=c'$.

\n Conversely, if (\ref{e1}) holds, then
$$\ln V_{F(\nu)}(m)=3\ln (1+\beta
m)+a\psi_{\nu}(m)+bk_{\nu}(\psi_{\nu}(m))+c.$$ Taking the
derivative, we deduce that the variance function satisfies the
differential equation
$$(1+\beta m) V'_{F(\nu)}(m)-3\beta V_{F(\nu)}(m)=(a+bm)(1+\beta m).$$ Solving
this equation by standard methods gives
\begin{eqnarray*}V_{F(\nu)}(m)&=&\lambda (1+\beta
m)^{3}+\frac{b}{\beta^{2}}(1+\beta m)^{2}-\displaystyle\frac{b-\beta
a}{2\beta^{2}}(1+\beta m),
\end{eqnarray*}
which is a polynomial of degree less than or equal to 3.
\end{Pff}

A third characterization of the cubic NEF's is based on a relation
between the associated families of prior distributions $\Pi^{\beta}$
and $\widetilde{\Pi}^{\beta}$.
\begin{theorem}\label{T2}Let $\nu$ be in $\mathcal{M}(\reel)$. $F(\nu)$ is cubic if and
 only if there exist $\beta$ in $B_{F(\nu)}$ such that
 $$k'_{\nu}(\Pi^{\beta})\subset
 \widetilde{\Pi}^{\beta}.$$\end{theorem}

\begin{Pff}
Suppose that $F(\nu)$ is cubic, then from Theorem \ref{T1} there
exist $\beta$ in $B_{F(\nu)}$ and $(a,b,c)$ in $\reel^{3}$ such
that
\begin{equation}\label{vrc}
V_{F(\nu)}(m)=(1+\beta \textit{m})^{3}\
\exp{(a\psi_{\nu}(m)+bk_{\nu}(\psi_{\nu}(m))+c)}.
\end{equation}
 Consider the set
$$\Omega=\{(t_{1},m_{1})\in\reel_{+}^{*}\times \ (M_{F(\nu)})_{\beta}\ ; t_{1}-b>0\ and\ \displaystyle\frac{t_{1}m_{1}+a}{t_{1}-b}\in
\ (M_{F(\nu)})_{\beta}\}.$$ We will show that
$k'_{\nu}(\Pi^{\beta})=\{\widetilde{\pi}^{\beta}_{t_{1},m_{1}};\
(t_{1},m_{1})\in\Omega\},$ which is a part of
$\widetilde{\Pi}^{\beta}$.\\Let $t>0$ and $m_{0}$ in
$(M_{F(\nu)})_{\beta}$, and denote by $\sigma$ the image by
$k'_{\nu}$ of the prior $\pi^{\beta}_{t,m_{0}}$ on $\theta$ defined
in (\ref{d}). We easily verify that
$$\sigma(\textit{dm})=C^{\beta}_{t,m_{0}}\ (1+\beta m)\ V_{F(\nu)}^{-1}(m)\
\exp( tm_{0}\psi_{\nu}(m)-tk_{\nu}(\psi_{\nu}(m)))\
\mathbf{1}_{(\textit{M}_{F(\nu)})_{\beta}}(\textit{m})\textit{dm}.$$
This using (\ref{vrc}) becomes
$$\sigma(\textit{dm})=\widetilde{C}^{\beta}_{t_{1},m_{1}}\ (1+\beta m)^{-2} \exp( t_{1}
m_{1}\psi_{\nu}(m)-t_{1}k_{\mu}(\psi_{\nu}(m)))\
\mathbf{1}_{(\textit{M}_{F(\nu)})_{\beta}}(\textit{m})\textit{dm},$$
with $ t_{1}=t+b$, $m_{1}=\displaystyle\frac{tm_{0}-a}{t+b}$ and
$\widetilde{C}^{\beta}_{t_{1},m_{1}}=C^{\beta}_{t_{1}-b,(t_{1}m_{1}+a)/(t_{1}-b)}.$\\We
have that $t_{1}-b=t>0$ and
$\displaystyle\frac{t_{1}m_{1}+a}{t_{1}-b}=m_{0}\in
(M_{F(\nu)})_{\beta},$ that is $(t_{1},m_{1})\in\Omega$. Hence
$$k'_{\nu}(\Pi^{\beta})\subset\{\widetilde{\pi}^{\beta}_{t_{1},m_{1}};\
(t_{1},m_{1})\in\Omega\}.$$ In the same, we verify that
$$\{\widetilde{\pi}^{\beta}_{t_{1},m_{1}};\
(t_{1},m_{1})\in\Omega\}\subset k'_{\nu}(\Pi^{\beta}).$$Finally, we
obtain that
$$k'_{\nu}(\Pi^{\beta})=\{\widetilde{\pi}^{\beta}_{t_{1},m_{1}};\
(t_{1},m_{1})\in\Omega\}\subset \widetilde{\Pi}^{\beta}$$
Conversely, suppose that $k'_{\nu}(\Pi^{\beta})\subset
\widetilde{\Pi}^{\beta}.$ The image of an element
$\pi^{\beta}_{t,m_{0}}$ of $\Pi^{\beta}$ by $k'_{\nu}$ is by the
very definition
$$k'_{\nu}(\pi^{\beta}_{t,m_{0}})(\textit{dm})=C^{\beta}_{t,m_{0}}\ (1+\beta m)\ (V_{F(\nu)}(m))^{-1}\
\exp (tm_{0}\psi_{\nu}(m)-tk_{\nu}(\psi_{\nu}(m)))\
\mathbf{1}_{(\textit{M}_{F(\nu)})_{\beta}}(\textit{m})\textit{dm}.$$
Since it is assumed to be in $\widetilde{\Pi}^{\beta}$, there exists
$(t_{1},m_{1})$ in $\reel^{*}_{+}\times (M_{F(\nu)})_{\beta}$ such
that
$$k'_{\nu}(\pi^{\beta}_{t,m_{0}})(\textit{dm})=\widetilde{C}^{\beta}_{t_{1},m_{1}}\
(1+\beta m)^{-2}\ \exp
(t_{1}m_{1}\psi_{\nu}(m)-t_{1}k_{\nu}(\psi_{\nu}(m)))\
\mathbf{1}_{(\textit{M}_{F})_{\beta}}(\textit{m})\textit{dm}.$$
Comparing these two expressions of
$k'_{\nu}(\pi^{\beta}_{t,m_{0}})$ gives
$$ V_{F(\nu)}(m)=(1+\beta m)^{3}\
\exp(a\psi_{\nu}(m)+bk_{\nu}(\psi_{\nu}(m))+c),$$ where
$$a=tm_{0}-t_{1}m_{1},\ b=t_{1}-t,\ and\ c=\ln
\left(\displaystyle\frac{C^{\beta}_{t,m_{0}}}{\widetilde{C}^{\beta}_{t_{1},m_{1}}}\right),$$
According to Theorem \ref{T1}, this is the desired result and the
proof is complete.
\end{Pff}

\section{Example}
In this section we illustrate our results by an example involving
the most famous family with variance function of degree 3 which is
the inverse Gaussian natural exponential family. \n Consider the
distribution
$$\nu(dx)=\displaystyle\frac{(1+x)^{-3/2}}{\sqrt{2\pi}}\
\exp(-\frac{x^{2}}{2(1+x)})\
\mathbf{1}_{]-1,+\infty[}(\textit{x})\textit{dx},$$ which is up to
an affine transformation an inverse Gaussian distribution. The NEF
generated by $\nu$ is given by
\begin{equation}\label{f}F(\nu)=\{\exp((1+x)\theta+\sqrt{-2\theta})\
\nu(dx);\ \theta<0\},\end{equation}its mean parametrization is
\begin{equation}\label{fm}F(\nu)=\{\exp(-\frac{1+x}{2(1+m)^{2}}+\frac{1}{1+m})\
\nu(dx);\ m>-1\}.\end{equation}For all $m>-1$ the variance
function is given by
\begin{equation}\label{vp}V_{F(\nu)}(m)=
(1+m)^{3}.\end{equation} Let $\mu=T_{-\beta}(\nu)$, we have
$$\mu(dx)=\frac{1}{\sqrt{2\pi}\sqrt{1-\beta x+x}}\
 \exp(-\frac{x^{2}}{2(1+(1-\beta) x)})\ \mathbf{1}_{\{1+(1-\beta) x)>0\}}(\textit{x})\textit{dx}$$
 Now for $\beta\neq 0$ we
 have
$$P(\beta,\theta,\nu)(dx)=\frac{e^{\theta+\sqrt{-2\theta}}}{\sqrt{2\pi}\sqrt{1-\beta
x+x}}\
\exp((\theta+\beta(-\theta-\sqrt{-2\theta}))x-\frac{x^{2}}{2(1-\beta
x+x)})\ \mathbf{1}_{\{\textit{1}-\beta
\textit{x}+\textit{x}>\textit{0}\}}(\textit{x})\textit{dx},$$ so
that
$$F^{\beta}=T_{-\beta}(F(\nu))=\{P(\beta,\theta,\nu)(dx);\
\theta<0\},$$ The corresponding family $\Pi^{\beta}$ of conjugate
prior distributions is the family of distributions
$$\pi^{\beta}_{t,m_{0}}(d\theta)=C^{\beta}_{t,m_{0}}\
(1+\beta(\frac{1}{\sqrt{-2\theta}}-1))\
\exp\{tm_{0}\theta+t(\theta+\sqrt{-2\theta})\}\
\mathbf{1}_{(\Theta)_{\beta}}(\theta)\textit{d}\theta,$$where
\begin{eqnarray*}(\Theta)_{\beta}&=&]-\frac{1}{2}(\frac{\beta}{\beta-1})^{2}\
,0[,\\
(\textit{M}_{F(\nu)})_{\beta}&=&\left\{\begin{array}{ccc}
 ]-1/\beta,+\infty[\ \ &if&\
\beta<0\\]-1,+\infty[\ \ &if&\
\beta=0\\]inf(-1,-1/\beta),+\infty[\ \ &if&\ \beta>0\
\\\end{array}\ \ \right.\end{eqnarray*} and
$$(C^{\beta}_{t,m_{0}})^{-1}=-\frac{\beta}{t}\
[1-\exp(-\frac{tm_{0}}{2}(\frac{\beta}{\beta-1})^{2}+\frac{t\beta}{\beta-1}-\frac{t}{2}(\frac{\beta}{\beta-1})^{2})]-(1+\beta
m_{0})\
\frac{1}{\frac{t}{\beta}-\frac{tm_{0}}{2}(\frac{\beta}{\beta-1})^{2}}.$$
defined for $t>0$ and $m_{0}$ in $(M_{F(\nu)})_{\beta}$\\
Also, in this example, the family $\widetilde{\Pi}^{\beta}$ is the
set of distributions defined for $t_{1}>0$ and $m_{1}$ in
$(M_{F(\nu)})_{\beta}$ by
$$\widetilde{\pi}^{\beta}_{t_{1},m_{1}}(dm)=\widetilde{C}^{\beta}_{t_{1},m_{1}}\
(1+\beta m)^{-2}\
\exp(-\frac{t_{1}(m_{1}+1)}{2(1+m)^{2}}+\frac{t_{1}}{1+m})\
\mathbf{1}_{(\textit{M}_{F(\nu)})_{\beta}}(\textit{m})\textit{dm},$$
 To see how Theorem \ref{T1} holds in this example, we need only to take
$\beta=1.$ Then $\mu=T_{-1}(\nu)$ is the standard gaussian
distribution with $V_{F(\mu)}(m')=1$ for $m'\in \reel$. We see that
$$V_{F(\nu)}(m)=(1+m)^{3}\
\exp(a\psi_{\nu}(m)+bk_{\nu}(\psi_{\nu}(m))+c),$$ with $a=b=c=0.$
In fact
$$V_{F(\mu)}(m')=\exp(a'\psi_{\mu}(m')+b'k_{\mu}(\psi_{\mu}(m'))+c'),$$
with $a'=b'=c'=0$, and using the relations $a=a',\ b=b'+\beta a',\
c=c',$ we get $a=b=c=0.$\\Concerning Theorem \ref{T2}, we fist
observe that the hypotheses in this theorem are well verified. In
fact, let $\pi^{1}_{t,m_{0}}$ be the prior on the natural parameter
$\theta$ and $\widetilde{\pi}^{1}_{t_{1},m_{1}}$ be the prior
on the mean parameter $m$, then\begin{eqnarray*}t_{1}&=&t+b=t>0\\
m_{1}&=&\frac{tm_{0}-a}{t+b}=m_{0}\in\
(M_{F(\nu)})_{1}=]1,+\infty[.\end{eqnarray*} The density function of
$\pi^{1}_{t,m_{0}}$ is equal to
$$\frac{1}{\sqrt{-2\theta}}\ \exp\{tm_{0}\theta+t(\theta+\sqrt{-2\theta})\},$$and for all $m>-1$ the density function of
$k'_{\nu}(\pi^{1}_{t,m_{0}})$ is given by
$$(1+m)^{-2}\
\exp\{-\frac{tm_{0}}{2(1+m)^{2}}-t(\frac{1}{2(1+m)^{2}}-\frac{1}{1+m})\}$$
which is equal to $\widetilde{\pi}^{1}_{t,m_{0}}.$

\end{document}